\documentclass{article}
\usepackage[english]{babel}
\usepackage{amsmath}
\usepackage{amsthm}
\usepackage{amssymb}
\usepackage{amscd}
\usepackage[latin1]{inputenc}

\newtheorem{teor}{Theorem}[section]

\newtheorem{lema}[teor]{Lemma}
\newtheorem{prop}[teor]{Proposition}
\newtheorem{cor}[teor]{Corollary}
\newtheorem{rem}[teor]{Remark}
\newtheorem{ejem}[teor]{Example}

\newtheorem{quesrem}[teor]{Question and Remark}

\begin{document}

\title{A duality theorem for  generalized Koszul algebras}

\author{Roberto Martínez Villa\\
Instituto de Matemáticas, UNAM, AP 61-3 \\
58089 Morelia, Michoacán\\ MEXICO \\ {\it mvilla@matmor.unam.mx} \vspace*{0.5cm} \and Manuel Saorín\\ Departamento de Matemáticas\\
Universidad de Murcia, Aptdo. 4021\\
30100 Espinardo, Murcia\\
SPAIN\\ {\it msaorinc@um.es}}

\date{}

\thanks{
The first named author thanks CONACYT for funding the research
project. The second one thanks the D.G.I. of the Spanish Ministry
of Education and Science and the Fundación "Séneca" of Murcia for
their financial support}

\maketitle

\begin{abstract}

{\bf
 We show that if $\Lambda$ is a $n$-Koszul algebra and $E=E(\Lambda )$ is its Yoneda algebra,
 then there is a full subcategory $\mathcal{L}_E$ of the category $Gr_E$ of graded $E$-modules,  which
 contains all the graded $E$-modules presented in even degrees,  that embeds fully faithfully in the
 category $C(Gr_\Lambda )$ of cochain complexes of graded $\Lambda$-modules. That extends the known equivalence,
 for $\Lambda$ Koszul (i.e. $n=2$), between $Gr_E$ and the category of linear complexes of  graded $\Lambda$-modules}

\end{abstract}

\section{Introduction}

From the classification of coherent sheaves over projective spaces
by Bernstein, Gelfand and Gelfand (see \cite{BG} or \cite{GM}),
Koszul algebras have deserved a lot of attention. A systematic
treatment of them was given in \cite{BGS}, where the authors
showed the existence of an equivalence of categories between large
subcategories of the graded derived categories of a Koszul algebra
$\Lambda$ and its (quadratic) dual $\Lambda^!$. They showed in
addition that $\Lambda^!$ is also Koszul and canonically
isomorphic to the Yoneda algebra $E=E(\Lambda )$ of $\Lambda$.
Recently (cf. \cite{MVS}),  the authors of the present paper
showed that there is an abelian version of the mentioned
equivalences, valid for more general graded algebras. Namely,
after defining $A^!$ for an arbitrary positively graded algebra
$A$, we showed that there is an equivalence between the category
$Gr_{A!}$ of graded $A!$-modules and the category $\mathcal{L}C_A$
(resp. $\mathcal{L}C^*_A$) of
 linear complexes of projective (resp. almost injective) graded
$A$-modules. In case $A$ is Koszul, the equivalence of \cite{BGS}
can be obtained from that by derivation.

As a generalization of Koszul algebras, Berger (\cite{B})
introduced $n$-Koszul algebras, where $n>1$ is an integer. They
form a class of $n$-homogeneous algebras which includes all
Artin-Schelter regular algebras of global dimension $\leq 3$ and,
in case $n=2$, coincides with the class of Koszul algebras.  If
$\Lambda$ is $n$-Koszul, then one also has a ($n$-homogeneous)
dual algebra, still denoted $\Lambda^!$, from which the Yoneda
algebra $E=E(\Lambda )$ can be obtained by killing supports and
appropriate regrading on $\Lambda^!$ (see \cite{GMMVZ}). However,
the existence, in the flavour of \cite{MVS}, of an equivalence of
categories between reasonably large subcategories of $Gr_E$ and of
the category $C(Gr_\Lambda )$ of cochain complexes of graded
$\Lambda$-modules is still unknown. The same lack of knowledge can
be applied to the existence, in the flavour of \cite{BGS},  of
equivalences between reasonably large subcategories of the derived
categories of graded modules over $E$ and $\Lambda$. Such an
equivalence exists, however, in the context of
$\mathbf{A}_\propto$-algebras (cf. \cite{K}).

The aim of this paper is to present, for every $n$-Koszul algebra
$\Lambda$ and each integer $m$, an equivalence between a full
subcategory $\mathcal{L}_E$ of $Gr_E$ containing all the graded
$E$-modules presented in even degrees and a full subcategory
$\mathcal{Y}$ (depending on $m$) of $C(Gr_\Lambda )$ (cf.
Corollary \ref{modules-over-Yoneda}). Since
$E=\Lambda^!_\mathcal{U}$ as ungraded algebras, where
$\mathcal{U}=n\mathbf{Z}\cup (n\mathbf{Z}+1)$, our strategy
consists in showing an equivalence of categories between suitable
subcategories of $Gr_{\Lambda^!_\mathcal{U}}$ and $C(Gr_\Lambda)$
(cf. Theorem \ref{grmodules-to-2complexes}), which is valid for
any graded algebra with relations of degree $\geq n$ and from
which the desired result follows easily. The proof of last theorem
is based on some consequences of the results in \cite{MVS} (see
Section 2), on  equivalences between subcategories of
$Gr_{\Lambda^!}$ and $Gr_{\Lambda^!_\mathcal{U}}$ which were given
in \cite{MVS2} for arbitrary group-graded algebras, and on the
transport of a torsion theory from $Gr_{\Lambda^!}$ to the
category $_n\mathcal{LC}_\Lambda$ of linear $n$-complexes of
almost injective graded $\Lambda$-modules (see Section 3).

In this paper we borrow the terminology from \cite{MVS}, with some
concrete adaptations. So the term {\bf positively graded algebra}
will stand for a $\mathbf{Z}$-graded algebra
$A=\oplus_{i\in\mathbf{Z}}A_i$ such that $A_i=0$, for all $i<0$,
$A_0$ is isomorphic to a direct product of copies of the ground
field $K$
 and $dim_KA_1<\propto$. The category of its right (resp. left)
 graded modules is denoted by $Gr_A$ (resp. $_AGr$) and the full
 subcategory of locally finite graded modules is
 denoted $lfgr_A$ (resp. $_Alfgr$). Unless explicitly said
 otherwise, 'module' will mean 'right module'. Notice that, in our
 situation,  the graded Jacobson radical of $A$ is $J^{gr}(A)=\oplus_{i>0}A_i$
and
 an
 object $T\in Gr_A$ is semisimple if, and only if,
 $TJ^{gr}(A)=0$. That allows, in particular, to identify the
 category $Mod_{A_0}$ with that of semisimple graded $A$-modules concentrated
 in degree $0$. If $M\in Gr_A$ then the
 {\bf graded socle} $Soc^{gr}(M)$, which is by definition the largest graded
 semisimple submodule of $M$,  is given by $Soc^{gr}(M)=\{x\in M:$
 $xJ^{gr}(A)=0\}$. If $\mathcal{X}\subseteq\mathbf{Z}$ then we
 shall say that $M$ is {\bf generated (resp. cogenerated) in degrees belonging
 to} $\mathcal{X}$ in case every nonzero factor (resp. subobject)
 of $M$ in $Gr_A$ has a support which intersects $\mathcal{X}$
 nontrivially. Also, we shall say that $M$ is {\bf presented in degrees belonging
 to} $\mathcal{X}$ in case it is the cokernel  of a
 morphism in $Gr_A$ between projective objects generated in
 degrees belonging to $\mathcal{X}$.
Notice that if $j\in\mathbf{Z}$ and
 $M$ is cogenerated in degree $j$, then   $Soc^{gr}(M)=M_j$.  Recall also that the {\bf projective graded modules} (i.e.
 the projective objects of $Gr_A$) are those in
 $Add(\oplus_{k\in\mathbf{Z}}A[k])$, while the {\bf almost injective
 graded
 modules} are those in $Add(\oplus_{k\in\mathbf{Z}}D(A)[k])$,
 where $D=\mathcal{H}om_{A_0}(-,A_0):_AGr\longrightarrow Gr_A$ is
 the canonical contravariant functor, which induces by restriction a duality
 $_Alfgr\stackrel{\cong^o}{\longrightarrow}lfgr_A$ (see \cite{MVS}
 for  details).

 We have a  finite
 quiver $Q$ associated to $A$, uniquely determined by the existence of isomorphisms $KQ_0\cong A_0$ and $KQ_1\cong A_1$
 (of $K$-algebras and $KQ_0-KQ_0-$bimodules, respectively),  together with a structural
 homomorphism of graded algebras $\pi_A :KQ\longrightarrow A$ whose
 image is the (graded) subalgebra $\bar{A}$ of $A$ generated by $A_0\oplus
 A_1$. When $\bar{A}=A$ (i.e. $\pi_A$ is surjective), we say that
 $A$ is a {\bf graded factor of a path algebra}. We shall use the
 letter $\Lambda$ instead of $A$ when we want to emphasize that
 the algebra is a graded factor of a path algebras. All tensors
 $\otimes$ in the paper are tensor over $KQ_0=A_0$.
In general, if $I=Ker(\pi_A)$ then  $I_k=\{x\in I:$ $x\text{ is
homogeneous of degree
 }k\}$ is an $A_0-A_0-$sub-bimodule of $KQ_k$, for every $k\geq 0$.  We can consider
 the orthogonal $I_n^\perp\subseteq KQ_n^{op}$ (for our fixed $n\geq 2$) with respect to the
 canonical duality $KQ_n^{op}\otimes KQ_n\longrightarrow KQ_0=A_0$
 (see \cite{MVS}). The graded algebra
$A^!=KQ^{op}/<I_n^\perp>$ will be called the {\bf n-homogeneous
dual algebra} of $A$ and we also put $^!A=(A^!)^{op}$. If $\Lambda
=KQ/I$ is a graded factor of a path algebra, we shall say that it
is {\bf n-homogeneous} (resp. {\bf has relations of degree}
$\geq${\bf n}) when $I$ is generated by $I_n$ (resp.
$\coprod_{k\geq n}I_k$).

\section{Graded modules versus n-complexes}

In this section we extend the equivalences of \cite{MVS} from
($2$-)complexes to $n$-complexes. Recall that we have a canonical
$\mathbf{Z}\times\mathbf{Z}$-grading on $A[X]$ by putting
$A[X]_{(i,j)}=A_iX^j$, in case $i,j\geq 0$, and $A[X]_{(i,j)}=0$
otherwise. We refer to that paper to see the interpretation of the
objects of  $Gr_{A[X]}$ as pairs $(P^\cdot ,d^\cdot )$, where
$P^\cdot\in Gr_A^\mathbf{Z}$ and $d^\cdot :P^\cdot\longrightarrow
P^\cdot$ is a morphism in $Gr_A^\mathbf{Z}$ of degree $+1$. We
then denote by $_n\mathcal{LC}_A$ (resp. $_n\mathcal{LC}_A^*$) the
full subcategory of $Gr_{A[X]}$ consisting of those pair $(P^\cdot
,d^\cdot )$ (resp. $(I^\cdot ,d^\cdot )$) satisfying the following
two conditions:

\begin{enumerate}
\item $(P^\cdot ,d^\cdot )$ (resp. $(I^\cdot ,d^\cdot )$) is a
$n$-complex (i.e. $d^{k+n-1}\circ ...\circ d^{k+1}\circ d^k=0$,
for all $k\in\mathbf{Z}$) \item
 $P^k$ (resp. $I^k$) is projective and generated in
degree $-k$ (resp. almost injective and cogenerated in degree
$-k$), for every $k\in\mathbf{Z}$
\end{enumerate}

The objects of $_n\mathcal{LC}_A$ (resp. $_n\mathcal{LC}_A^*$)
will be called {\bf linear n-complexes of projective (resp. almost
injective)} graded $A$-modules. We shall denote by
$_n\mathcal{L}c_A$ (resp. $_n\mathcal{L}c_A^*$) the full
subcategory of $_n\mathcal{LC}_A$ (resp. $_n\mathcal{LC}_A^*$)
consisting of those $(P^\cdot ,d^\cdot )$ (resp. $(I^\cdot
,d^\cdot )$) such that $P^k$ (resp. $I^k$) is finitely generated
(resp. finitely cogenerated), for all $k\in\mathbf{Z}$. We then
consider the fully faithful embddings $\Psi_A,
\nu_A:_{KQ}Gr\longrightarrow Gr_{A[X]}$ given in
\cite{MVS}[Theorems 2.4 and 2.10]. We refer to that paper for
their explicit definition, which we shall freely use here.

\begin{prop} \label{gradedmodules-ncomplexes1}
Let $A=\oplus_{i\geq 0}A_i$ be a positively graded algebra with
quiver $Q$ and $A^!$ be its $n$-homogeneous dual. Then  $\Psi
=\Psi_A :_{KQ}Gr\longrightarrow Gr_{A[X]}$ induces by restriction
equivalences
$_{^!A}Gr=Gr_{A^!}\stackrel{\cong}{\longrightarrow}_n\mathcal{LC}_A$
and $_{^!A}lfgr
=lfgr_{A^!}\stackrel{\cong}{\longrightarrow}_n\mathcal{L}c_A$.
\end{prop}

\begin{proof}
By the proof of \cite{MVS}[Theorem 2.4], we know that $\Psi$
establishes an equivalence of categories
$_{KQ}Gr\stackrel{\cong}{\longrightarrow}\mathcal{LG}_A$, where
$\mathcal{LG}_A$ is the full subcategory of $Gr_{A[X]}$ consisting
of those objects $(P^\cdot ,d^\cdot )$ such that $P^k$ is a
projective object of $Gr_A$ generated in degree $-k$, for all
$k\in\mathbf{Z}$. We only need to prove that if $M\in
_{KQ}Gr=Gr_{KQ^{op}}$ and $\Psi (M)=(P^\cdot ,d^\cdot )$, then
$M\cdot I_n^\perp =0$ if, and only if, $(P^\cdot ,d^\cdot )$ is a
$n$-complex.

We consider  the canonical $K$-algebra homomorphism $\pi_A:KQ
\longrightarrow A$ whose kernel is $I$.  It is convenient in the
rest of this proof to view graded left $KQ$-modules as graded
right $KQ^{op}$-modules. On one hand, $\psi (M)$ is a cochain
$n$-complex iff the composition $M_k\otimes
A[k]\stackrel{d^k}{\longrightarrow} M_{k+1}\otimes
A[k+1]...\stackrel{d^{k+n-1}}{\longrightarrow}M_{k+n}\otimes
A[k+n]$ is zero, for each fixed $k\in\mathbf{Z}$. But that is
equivalent to say that $d^{k+n-1}\circ...\circ d^k$ vanish on
$M_k\cong M_k\otimes A_0$.
 Direct
calculation shows that $(d^{k+n-1}\circ...\circ d^k)(x)=
\sum_{p\in Q_{n}} xp^o\otimes \bar{p}$,  for all $x\in M_k$, where
$\bar{p}=\pi_A(p)$. Our goal is to show that this latter sum is
zero for all $x\in M_k$ iff $M_k\cdot I_{n}^\perp=0$. To do that
we consider an ordering of
$Q_{n}=\{p_1,...,p_r,...p_{r+s},...p_{r+s+t}\}$, where i)
$\{p_1,...,p_r\}$ is a basis of  $KQ_{n}$ modulo $I$; ii)
$\{p_{r+1},...p_{r+s}\}$  gathers the remaining $p\in Q_{n}$ which
do not belong to $I$; iii) $\{p_{r+s+1},...p_{r+s+t}\}$ gathers
the $p\in Q_{n}$ which belong to $I$. Since $0\neq\pi_A(p_j)\in
A_1^{n}\cong\frac{KQ_n}{I_n}$ when $j=r+1,...,r+s$, we have
uniquely determined linear combinations $p_j-\sum_{1\leq i\leq
r}\lambda_{ij}p_i$ belonging to $I$, for $j=r+1,...,r+s$, which
can be taken with the property that $\lambda_{ij}\neq 0$ implies
that $p_i$ and $p_j$ share origin and terminus.
  Those linear combinations together with the  $p_j$,
with $j=r+s+1,...,r+s+t$,  form a basis of $I_{n}$. By canonical
methods of Linear Algebra, a basis of $I_{n}^\perp$ is then given
by the elements of the form
$h_i=p_i^o+\lambda_{i,r+1}p_{r+1}^o+...+\lambda_{i,r+s}p_{r+s}^o$
($i=1,...,r$). Now we have  $(d^{k+n-1}\circ ...\circ d^k)(x)=
\sum_{1\leq i\leq r+s+t}xp_i^o\otimes\bar{p}_i=$ $\sum_{1\leq
i\leq r} xp^o_i\otimes\bar{p}_i$ + $\sum_{r+1\leq j\leq
r+s}xp^o_j\otimes(\sum_{1\leq i\leq r}\lambda_{ij} \bar{p}_i)$
(notice that, for $r+s+1\leq i\leq r+s+t$,  the summand
$xp^o_i\otimes\bar{p}_i$ is zero because $p_i\in I$). We can write
the last summatory as $\sum_{1\leq i\leq r} x(p^o_i+\sum_{r+1\leq
j\leq r+s}\lambda_{ij}p^o_j)\otimes\bar{p}_i= \sum_{1\leq i\leq r}
xh_i\otimes \bar{p}_i$. The fact that $\{\bar{p}_i:$ $i=1,...,r\}$
is  a $K$-linearly independent subset of $A_{n}$ easily implies
that the last summatory is zero in $ M_{k+n}\otimes A_{n}$ iff
$xh_i=0$ for $i=1,...,r$. This is equivalent to say that $M_k\cdot
I_{n}^\perp=0$ and we are done.
\end{proof}

If the above result completes \cite{MVS}[Theorem 2.4], the
following one completes \cite{MVS}[Theorem 2.10]:

\begin{prop} \label{gradedmodules-ncomplexes2}
Let $A=\oplus_{i\geq 0}A_i$ be a locally finite positively graded
algebra with quiver $Q$ and $A^!$ be its $n$-homogeneous dual.
Then  $\upsilon =\upsilon_A:_{KQ}Gr\longrightarrow Gr_{A[X]}$
induces by restriction equivalences
$_{^!A}Gr=Gr_{A^!}\stackrel{\cong}{\longrightarrow}_{n}\mathcal{LC}_A^*$
and  $_{^!A}lfgr
=lfgr_{A^!}\stackrel{\cong}{\longrightarrow}_{n}\mathcal{L}c_A^*$.
\end{prop}

\begin{proof}

Let  $\pi :KQ\longrightarrow A$ be the canonical homomorphism of
graded algebras whose kernel is $I$. Again, we view the graded
left $KQ$-modules as right $KQ^{op}$-modules. We need to prove
that  if $N\in Gr_{KQ^{op}}$ then $\upsilon (N)$ is a cochain
$n$-complex iff $N_k\cdot I_{n}^\perp=0$, for all
$k\in\mathbf{Z}$.  We have that $\upsilon (N)$ is a cochain
$n$-complex iff the composition
$\mathcal{H}om_{A_0}(A,N_k)[k]\stackrel{d^k}{\longrightarrow}
\mathcal{H}om_{A_0}(A,N_{k+1})[k+1]...$
\newline $\stackrel{d^{k+n-1}}{\longrightarrow}\mathcal{H}om_{A_0}(A,N_{k+n})[k+n]$
is zero, for all $k\in\mathbf{Z}$, iff its $-(k+n)$-component
$Hom_{A_0}(A_{n},N_k) \stackrel{d}{\longrightarrow}
Hom_{A_0}(A_{n-1},N_{k+1})\stackrel{d}{\longrightarrow}...
\stackrel{d}{\longrightarrow}\mathcal{H}om_{A_0}(A_0,N_{k+n})\cong
N_{k+n}$ is zero, for all $k\in\mathbf{Z}$. The latter happens iff
$\sum_{p\in Q_{n}}f(\bar{p})p^o =0$, for all $f\in
Hom_{A_0}(A_{n},N_k)$, where $\bar{p}=\pi (p)$. Since
$A_1^{n}=A_1\cdot A_1...\cdot A_1$ is a direct summand of $A_{n}$
in $Mod_{A_0}$ and $\bar{p}\in A_1^{n}$, for all $p\in Q_{n}$, we
get that $\upsilon (N)$ is a cochain $n$-complex iff $\sum_{p\in
Q_{n}}f(\bar{p})p^o=0$, for all $f\in Hom_{A_0}(A_1^{n},N_k)$. Now
we choose an ordering of
$Q_{n}=\{p_1,...,p_r,...p_{r+s},...p_{r+s+t}\}$ with the same
criterion as in the proof of Proposition
\ref{gradedmodules-ncomplexes1}. Then
$\{\bar{p}_1,...,\bar{p}_r\}$ is a basis of $A_1^{n}$ and we have
$\sum_{p\in Q_{n}}f(\bar{p})p^o= \sum_{1\leq i\leq
r}f(\bar{p}_i)p^o_i+ \sum_{r+1\leq j\leq r+s}f(\sum_{1\leq i\leq
r}\lambda_{ij}\bar{p}_i)p^o_j= \sum_{1\leq i\leq
r}f(\bar{p}_i)(p^o_i+ \sum_{r+1\leq j\leq r+s}\lambda_{ij}p^o_j)
=\sum_{1\leq i\leq r}f(\bar{p}_i)h_i$, for every \newline $f\in
Hom_{A_0}(A_1^{n},N_k)$, with the same terminology of the proof of
Proposition \ref{gradedmodules-ncomplexes1}. On the other hand,
the $A_0$-homorphisms $A_1^{n}\longrightarrow N_k$ of the form
$x\bar{p}_i^{*}(-):a\longrightarrow x\bar{p}_i^{*}(a)$, with
$i=1,...,r$ and $x\in N_k$,  generate $Hom_{A_0}(A_1^{n},N_k)$
(see \cite{MVS}[Remark 2.1]). But, for $f=x\bar{p}_s^{*}(-)$, we
have $\sum_{1\leq i\leq r}f(\bar{p}_i)h_i=xh_s$. Consequently,
$\upsilon (N)$ is a $n$-complex iff $xh_s=0$, for all $s=1,..,r$
and $x\in N_k$, that is, iff $N_k\cdot I_{n}^\perp=0$ for all
$k\in\mathbf{Z}$. That ends the proof.
 \end{proof}

The reader is invited to extend to general $n$-homogeneous
algebras results like the equivalence of assertions 1, 2 and 5 in
\cite{MVS}[Corollary 3.4] and of assertions 1 and 2 of
\cite{MVS}[Corollary 3,5], which were  given there for quadratic
algebras.

\section{Transport of a torsion theory}
We know from \cite{MVS2} that if $A=\oplus_{i\in\mathbf{Z}}A_i$ is
a $\mathbf{Z}$-graded algebra and $\mathcal{S}\subseteq\mathbf{Z}$
is any subset, then $\mathcal{T}_\mathcal{S}=\{M\in Gr_A:$
$M_\mathcal{S}=0\}$ is a hereditary torsion class  closed for
products in $Gr_A$. The following is a handy way of identifying
its associated torsionfree class.

\begin{lema} \label{torsionfree}
Suppose that $A=\oplus_{i\geq 0}A_i$ is positively graded,
generated in degrees 0,1 and $A_0$ is semisimple. If
$\mathcal{S}\subseteq\mathbf{Z}$ is not upper bounded then, for a
graded $A$-module $M=\underset{n\in \mathbf{Z}}{\oplus }M_{n}$,
the following statements are equivalent:

\begin{enumerate}
\item $M$ is $\mathcal{T}$- torsionfree, where
$\mathcal{T}=\mathcal{T}_\mathcal{S}$  \item
$Hom_{Gr_A}(A_0[-m],M)=0$, for all
$m\in\mathbf{Z}\setminus\mathcal{S}$
\end{enumerate}
\end{lema}

\begin{proof}
$1)\Longrightarrow 2)$ is clear

$2)\Longrightarrow 1)$ Notice that every morphism
$A_0[-m]\longrightarrow M$ in $Gr_A$ is given by left
multiplication by an element $x\in ann_{M_m}(A_1)$. Hence the
hypothesis is equivalent to say that $ann_{M_m}(A_1)=0$, for all
$m\in\mathbf{Z}\setminus\mathcal{S}$.  Recall that
$Supp(t(M))\subseteq\mathbf{Z}\setminus\mathcal{S}$. If $t(M)\neq
0$ ( equivalently, $Supp(t(M))\neq\emptyset$)  and $m\in
Supp(t(M))$, then  we pick up $x\in t(M)_{m}\smallsetminus \left\{
0\right\} .$ By assumption $xA_{1}\neq 0$ which implies that
$t(M)_{m+1}\neq 0$ and, hence, that $m+1\in Supp(t(M))$. By
recurrence we get that the interval $[m,+\propto )$ is contained
in $Supp(t(M))\subseteq\mathbf{Z}\setminus\mathcal{S}$. But that
contradicts the fact that $\mathcal{S}$ is not upper bounded.

\end{proof}

Throughout the rest of the  section $\Lambda =\oplus_{i\geq
0}\Lambda_i$ is a graded factor of a path algebra (see section 1).
We  want to transfer the results of \cite{MVS2}  (for
$A=\Lambda^!$) from $Gr_{\Lambda^!}$ to $_n\mathcal{L}C^*_\Lambda$
and $_n\mathcal{L}C_\Lambda$ via $\upsilon$ and $\Psi$. We refer
the reader to \cite{MVS2} for the definition and terminology about
 ring-supporting subsets and right modular pairs of subsets of a
group. Here we shall fix a modular pair
$(\mathcal{S},\mathcal{U})=(m+\mathcal{U},\mathcal{U})$, with
$\mathcal{U}= \bigcup_{k\in\mathbf{Z}}[kn,kn+r]$ and
$m\in\mathbf{Z}$, where $n\geq 2$ and $0\leq 2r\leq n$, with
strict inequality $2r<n$ in case $n>2$.
 Then, in case $2r=n=2$,  we have $\mathcal{S}=\mathcal{U}=\mathbf{Z}=(\mathcal{S}:\mathcal{U})$. In any other case,
  we have $(\mathcal{S}:\mathcal{U})=m+(\mathcal{U}:\mathcal{U})=m+n\mathbf{Z}$.

\begin{lema} \label{transfer-of-cogeneration}
Let $I^\cdot =(I^\cdot ,d^\cdot )\in _n\mathcal{L}C^*_\Lambda$ be
a linear $n$-complex of almost injective graded $\Lambda$-modules
and let $M\in Gr_{\Lambda^!}$ be such that $I^\cdot
=\upsilon_\Lambda (M)$. The following assertions hold:

\begin{enumerate}
 \item $M\in\mathcal{T}=\mathcal{T}_\mathcal{S}$ if, and only if,
 $I^j=0$ for all $j\in\mathcal{S}=\bigcup_{k\in\mathbf{Z}}[m+kn,m+kn+r]$  \item $M$ is
 $\mathcal{T}$-torsionfree  iff
 $Hom_{Gr_\Lambda}(D(\Lambda )[j], Ker(d^j))=0$ for all
 $j\not\in\mathcal{S}$ \item If $0\leq 2r<n$ then the following conditions are equivalent:

\begin{enumerate}
 \item[a)]$M$ is generated in degrees belonging to
$(\mathcal{S}:\mathcal{U})=m+n\mathbf{Z}$ \item[b)]
  The differencial $d^{j-1}:I^{j-1}\longrightarrow I^{j}$
satisfies that $(Im d^{j-1})_{-j}=(I^j)_{-j}$  for all
$j\not\equiv m$ (mod $n$) \item[c)] $Soc^{gr}(I^j)\subseteq
Im(d^{j-1})$, for all $j\not\equiv m$ (mod $n$).
\end{enumerate}
\end{enumerate}
\end{lema}

\begin{proof}
By definition of $\upsilon$, we have
$I^j=\mathcal{H}om_{\Lambda_0}(\Lambda,M_j)[j]$ and, hence,
assertion 1 follows.  The equivalence $\upsilon$ takes
$\Lambda_0[-j]$ onto the stalk $n$-complex  $D(\Lambda )[j]$ (at
the position $j$). Then, using Lemma \ref{torsionfree}, we have
that $M$ is $\mathcal{T}$-torsionfree if, and only if, there are
no nonzero morphisms $f^\cdot:D(\Lambda )[j]\longrightarrow
I^\cdot$ in $_n\mathcal{L}C^*_\Lambda$, for all
$j\not\in\mathcal{S}$. But such a morphism is completely
determined by the induced morphism in $Gr_\Lambda$ $f:D(\Lambda
)[j]\longrightarrow Ker(d^j)$. Therefore $M$ is
$\mathcal{T}$-torsionfree if, and only if,
$Hom_{Gr_\Lambda}(D(\Lambda )[j], Ker(d^j))=0$, for all
$j\not\in\mathcal{S}$, which proves assertion 2.

On the other hand, since the algebra $\Lambda^!$ is generated in
degrees $0,1$, it is easy to see that $M$ is generated in degrees
belonging to $(\mathcal{S}:\mathcal{U})=m+n\mathbf{Z}$ if, and
only if, the multiplication map
$M_{j-1}\otimes\Lambda^!_1\longrightarrow M_j$ is surjective, for
all $j\not\in m+n\mathbf{Z}$. By using adjunction,  that is
equivalent to say that the $-j$-component
$Hom_{\Lambda_0}(\Lambda_1,M_{j-1})\longrightarrow M_j$ of the
'differential' $d^{j-1}:I^{j-1}=\mathcal{H}om_{\Lambda_0}(\Lambda
,M_{j-1})[j-1]\longrightarrow\mathcal{H}om_{\Lambda_0}(\Lambda
,M_{j})[j]=I^j$ is surjective for all $j\not\equiv m$ (mod $n$).
Then the equivalence of conditions a) and b) in 3) follows. The
equivalence of b) and c) is clear since
$Soc^{gr}(I^j)=(I^j)_{-j}$, for all $j\in\mathbf{Z}$.
\end{proof}

Recall from \cite{MVS2} that
$\mathcal{G}(\mathcal{S},\mathcal{U})$ is the full subcategory of
$Gr_A$ with objects those $M\in Gr_A$ which are
$\mathcal{T}$-torsionfree and generated in degrees belonging to
$(\mathcal{S}:\mathcal{U})$. We have the following:

\begin{prop} \label{transfer-of-torsion-and-quotient-categories}
Let $\mathcal{T}=\mathcal{T}_\mathcal{S}$  be the hereditary
torsion class in $Gr_{\Lambda^!}$ defined by
$\mathcal{S}=\bigcup_{k\in\mathbf{Z}}[m+kn,m+kn+r]$. The following
assertions hold:

\begin{enumerate}
\item $\mathcal{T}^*=\upsilon_\Lambda (\mathcal{T})$ is the
hereditary torsion class of $_n\mathcal{L}C^*_\Lambda$ consisting
of those $I^\cdot\in _n\mathcal{LC}^*_\Lambda$ such that $I^j=0$,
for all $j\in\mathcal{S}$ \item The functor $\nu =\nu_\Lambda
:Gr_{\Lambda^!}\longrightarrow Gr_{\Lambda [X]}$ induces a
commutative diagram:

\vspace*{0.5cm}

\setlength{\unitlength}{1mm}
\begin{picture}(140,30)
\put(25,25){$\mathcal{G}(\mathcal{S},\mathcal{U})$}
\put(42,25){$\hookrightarrow$} \put(50,25){$Gr_{\Lambda^!}$}
\put(59,26){\vector(1,0){20}}
\put(81,25){$\frac{Gr_{\Lambda^!}}{\mathcal{T}}$}

\put(24,4){$\mathcal{G}^*(\mathcal{S},\mathcal{U})$}
\put(42,5){$\hookrightarrow$}
\put(49,4){$_n\mathcal{LC}_\Lambda^*$}
\put(59,5){\vector(1,0){20}}
\put(80,4){$\frac{_n\mathcal{LC}_\Lambda^*}{\mathcal{T^*}}$}

\put(30,23){\vector(0,-1){12}} \put(53,23){\vector(0,-1){12}}
\put(83,23){\vector(0,-1){12}}

\put(67,28){$pr.$}  \put(67,2){$pr.$}

\put(27,16){$\cong$} \put(50,16){$\cong$}
 \put(80,16){$\cong$}

\end{picture}

where the vertical arrows are equivalences of categories and  the
compositions of horizontal arrows are fully faithful embeddings.
Here
$\mathcal{G}^*(\mathcal{S},\mathcal{U})=_n\mathcal{LC}_\Lambda^*$,
when $2r=n=2$, and is the full subcategory of
$_n\mathcal{L}C^*_\Lambda$ consisting of those $I^\cdot\in
_n\mathcal{LC}^*_\Lambda$ satisfying  conditions a) and b) below,
when $0\leq 2r<n$:

\begin{enumerate}
\item $Hom_{Gr_\Lambda}(D(\Lambda)[j], Ker(d^j))=0$, for all
$j\not\in\mathcal{S}$ \item $Soc^{gr}(I^j)\subseteq Im(d^{j-1})$,
for all $j\not\equiv m$ (mod $n$)
\end{enumerate}

\end{enumerate}
\end{prop}

\begin{proof}
In case $2r=n=2$ we have
$\mathcal{G}(\mathcal{S},\mathcal{U})=Gr_A$ and the result follows
from \cite{MVS}[Theorem 2.10]. If $0\leq 2r<n$ then the result
follows from Proposition \ref{gradedmodules-ncomplexes2}, Lemma
\ref{transfer-of-cogeneration} and \cite{MVS2}[Theorem 2.7].
\end{proof}

\section{From $n$-complexes to ($2$-)complexes}

In this section we  consider the case $r=1$ of Section 3, i.e.,
the modular pair is
$(\mathcal{S},\mathcal{U})=(m+\mathcal{U},\mathcal{U})$, where
$\mathcal{U}=n\mathbf{Z}\cup(n\mathbf{Z}+1)$ and $m\in\mathbf{Z}$.
 We want to pass from $n$-complexes to
(2-)complexes via the appropriate contraction. We will consider
parallels of the canonical contraction (see, e.g., \cite{BDVW} and
\cite{BM}). Let $C(Gr_\Lambda )$ be the category of
($2$-)complexes of graded $\Lambda$-modules. The unique strictly
increasing function $\delta_m=:\delta_{(\mathcal{S},m)}
:\mathbf{Z}\longrightarrow\mathbf{Z}$ such that $\delta_{m}(0)=m$
and $Im(\delta_{m})=\mathcal{S}$ is given by $\delta_m (2k)=m+kn$
and $\delta_m(2k+1)=m+kn+1$, for all $k\in\mathbf{Z}$ (cf.
\cite{MVS2}[Lemma 4.9]). Hence, $\delta_{m}(j)=m+\delta (j)$ for
all $j\in\mathbf{Z}$, where $\delta=\delta_{(\mathcal{U},0)}$ is
the map used, for instance, in \cite{GMMVZ}. We have an obvious
additive functor
$H=H_m:_n\mathcal{L}\mathcal{C}^*_\Lambda\longrightarrow
C(Gr_\Lambda )$ defined as follows.  We take $H(I^\cdot
)=\tilde{I}^\cdot$ where $\tilde{I}^k=I^{\delta_m(k)}$, for all
$k\in\mathbf{Z}$ and, as differentials,
$\tilde{d}^{2j}=d^{m+jn}:\tilde{I}^{2j}=I^{m+jn}\longrightarrow
I^{m+jn+1}=\tilde{I}^{2j+1}$ and
$\tilde{d}^{2j+1}=d^{m+(j+1)n-1}\circ ...\circ
d^{m+jn+1}:\tilde{I}^{2j+1}=I^{m+jn+1}\longrightarrow
I^{m+(j+1)n}=\tilde{I}^{2j+2}$. The objects in the essential image
of $H$ will be called  $H${\bf -liftable}.

The following observation is trivial, but very useful.

\begin{rem} \label{remark}
Let $A=\oplus_{i\geq 0}A_i$ be a positively graded algebra
generated in degrees $0,1$ and $(\mathcal{S},\mathcal{U})$ as
above. If $M,N\in Gr_{A_\mathcal{U}}$ and
$f=(f_i:M_i\longrightarrow N_i)_{i\in\mathbf{Z}}$ is a family of
morphisms in $Mod_{A_0}$, then the following assertions are
equivalent:

\begin{enumerate}
\item $f$ is a morphism in $Gr_{A_\mathcal{U}}$ \item
$f(xa)=f(x)a$, for all $x\in M_i$ and $a\in A_1\cup A_n$, where
$i\in Supp(M)$
\end{enumerate}
\end{rem}

\begin{proof}
 It is a straightforward consequence of the fact that, as an algebra,
 $A_\mathcal{U}$ is generated by $A_0$, $A_1$ and $A_n$
\end{proof}

\begin{prop} \label{from-ncomplex-to-2complex1}
The functor
$H=H_m:_n\mathcal{L}\mathcal{C}^*_\Lambda\longrightarrow
C(Gr_\Lambda )$ is an exact functor having the following
properties:

\begin{enumerate}
\item $Ker(H)=\mathcal{T}^*_{\mathcal{S}}=:\mathcal{T}^*$ and $H$
induces a faithful functor
$\frac{_n\mathcal{L}\mathcal{C}_\Lambda}{\mathcal{T}^*}\longrightarrow
C(Gr_\Lambda)$ \item The composition
$\mathcal{G}^*(\mathcal{S},\mathcal{U})\hookrightarrow
_n\mathcal{L}C^*_\Lambda\stackrel{H}{\longrightarrow}C(Gr_\Lambda
)$ is a  faithful functor whose essential image consists of those
$H$-liftable  complexes $(\tilde{I}^\cdot ,\tilde{d}^\cdot )$
satisfying  condition a) below, in case $n=2$, or both conditions
a) and b), in case $n>2$:

\begin{enumerate}
\item $\tilde{I}^j$ is almost injective and cogenerated in degree
$-\delta_m(j)$, for all $j\in\mathbf{Z}$ \item
$Soc^{gr}(\tilde{I}^{2j+1})\subseteq Im(\tilde{d}^{2j})$, for all
$j\in\mathbf{Z}$
\end{enumerate}

\item When  $\Lambda$ has relations of  degree $\geq n$,  the
composition $\mathcal{G}^*(\mathcal{S},\mathcal{U})\hookrightarrow
_n\mathcal{L}C^*_\Lambda\stackrel{H}{\longrightarrow} C(Gr_\Lambda
)$ is also full
\end{enumerate}
\end{prop}

\begin{proof}
If $n=2$ then $\mathcal{T}^*=0$,
$\mathcal{G}^*(\mathcal{S},\mathcal{U})=\mathcal{LC}^*_\Lambda$
and $H$ is the composition $\mathcal{LC}^*_\Lambda\hookrightarrow
C(Gr_\Lambda )\stackrel{?[m]}{\longrightarrow}C(Gr_\Lambda )$,
where $?[m]$ is the canonical shifting of $(2-)$complexes. Then
the functor $H$ is a fully faithful embedding and the three
assertions trivially hold in this case.

We assume in the rest of the proof that $n>2$. The exactness of
$H$ and the fact that $Ker(H)=\mathcal{T}^*$ are clear. In order
to prove property 1, we consider  a morphism $f^\cdot
:I^\cdot\longrightarrow J^\cdot$ in
$_n\mathcal{L}\mathcal{C}^*_\Lambda$ and will prove that
$H(f^\cdot )=0$ iff $f^\cdot$ is the zero morphism in the quotient
category
$\frac{_n\mathcal{L}\mathcal{C}^*_\Lambda}{\mathcal{T}^*}$. To see
that, using the equivalence of categories $\upsilon=
\upsilon_\Lambda:Gr_{\Lambda^!}\stackrel{\cong}{\longrightarrow}
_n\mathcal{L}\mathcal{C}^*_\Lambda$ (cf. Proposition
\ref{gradedmodules-ncomplexes2}), we have uniquely determined
objects $M,N$ and morphism $g:M\longrightarrow N$ in
$Gr_{\Lambda^!}$ such that $I^\cdot =\upsilon (M)$, $J^\cdot
=\upsilon (N)$ and $f^\cdot =\upsilon (g)$. One readily sees that
$H(f^\cdot )=0$ iff the functor
$(-)_{\mathcal{S}}:Gr_{\Lambda^!}\longrightarrow Gr_K$ maps $g$
onto zero. But, by \cite{MVS2}[Proposition 2.1],  that happens iff
$g$ is the zero morphism in the quotient category
$\frac{Gr_{\Lambda^!}}{\mathcal{T}}$, where
$\mathcal{T}=\mathcal{T}_{\mathcal{S}}$. But, by Proposition
\ref{transfer-of-torsion-and-quotient-categories}, we have that
$g=0$ in $\frac{Gr_{\Lambda^!}}{\mathcal{T}}$ iff $f=0$ in
$\frac{_n\mathcal{L}\mathcal{C}^*_\Lambda}{\mathcal{T}^*}$.

On the other hand, also by Proposition
\ref{transfer-of-torsion-and-quotient-categories},  the
composition
$\mathcal{G}^*(\mathcal{S},\mathcal{U})\hookrightarrow$
$_n\mathcal{L}C^*_A\stackrel{pr}{\twoheadrightarrow}
\frac{_n\mathcal{L}C^*_\Lambda}{\mathcal{T}^*}$ is fully faithful.
That together with the above paragraph give that the composition
$\mathcal{G}^*(\mathcal{S},\mathcal{U})\hookrightarrow
_n\mathcal{L}C^*_\Lambda\stackrel{H}{\longrightarrow}C(Gr_\Lambda
)$ is a faithful functor. By definition of
$\mathcal{G}^*(\mathcal{S},\mathcal{U})$, we readily see that if
$(\tilde{I}\cdot ,\tilde{d}^\cdot )$ is the image of $(I^\cdot
,d^\cdot )\in\mathcal{G}^*(\mathcal{S},\mathcal{U})$  by $H$, then
$(Im(\tilde{d}^{2j}))_{-m-jn-1}=(\tilde{I}^{2j+1})_{-m-jn-1}$ for
all $j\in\mathbf{Z}$ (i.e. $Soc^{gr}(\tilde{I}^{2j+1})\subseteq
Im(\tilde{d}^{2j})$). Conversely, suppose that $(\tilde{I}\cdot
,\tilde{d}^\cdot )$ is an $H$-liftable (2-)complex satifying
conditions a) and b) of the statement of the proposition.
 We first choose $J^\cdot\in
_n\mathcal{L}\mathcal{C}_\Lambda$ such that $H(J^\cdot
)=\tilde{I}^\cdot$.  Then we have a unique $M\in Gr_{\Lambda^!}$
such that $\upsilon (M)=J^\cdot$. Now condition b) translates into
the fact that $(Im (d^{m+jn}))_{-m-jn-1}=(J^{m+jn+1})_{-m-jn-1}$,
for all $j\in\mathbf{Z}$ (here $d^\cdot$ is the 'differential' of
the $n$-complex $J^\cdot$). Bearing in mind \cite{MVS}[Lemma 2.9]
and the definition of $\nu$, that means that the multiplication
map $M_{m+jn}\otimes \Lambda^!_1\longrightarrow M_{m+jn+1}$ is
surjective, for all $j\in\mathbf{Z}$. From that we get that if
$M'$ is the graded $\Lambda^!$-submodule of $M$ generated by
$M_{m+n\mathbf{Z}}$, then
$Supp(\frac{M}{M'})\subseteq\mathbf{Z}\setminus\mathcal{S}$, so
that  $\frac{M}{M'}\in\mathcal{T}$ and, hence, $M\cong
M'\cong\frac{M'}{t(M')}$ in $Gr_{\Lambda^!}/\mathcal{T}$. Then
$J^\cdot =\upsilon (M)\cong \upsilon (\frac{M'}{t(M')})$ in
$\frac{_n\mathcal{L}C^*_\Lambda}{\mathcal{T}^*}$. We take $I^\cdot
=\upsilon (\frac{M'}{t(M')})$ and have that $H(I^\cdot
)=\tilde{I}^\cdot$. Since
$\frac{M'}{t(M')}\in\mathcal{G}(\mathcal{S},\mathcal{U})$, we
conclude that $I^\cdot\in\mathcal{G}^*(\mathcal{S},\mathcal{U})$
as desired.

 We next prove assertion 3  assuming that the relations for $\Lambda$
 have degree
 $\geq n$. Take
$I^\cdot ,J^\cdot\in\mathcal{G}^*(\mathcal{S},\mathcal{U})$. Then
we have uniquely determined $M,N\in
\mathcal{G}(\mathcal{S},\mathcal{U})$ such that $\upsilon
(M)=I^\cdot$ and $\upsilon(N)=J^\cdot$.  Let $f:H(I^\cdot
)\longrightarrow H(J^\cdot )$ be a morphism in $C(Gr_\Lambda )$.
Then we get morphisms $f^s:I^s=\mathcal{H}om_{\Lambda_0}(\Lambda
,M_s)[s]\longrightarrow \mathcal{H}om_{\Lambda_0}(\Lambda
,N_s)[s]=J^s$ in $Gr_\Lambda$, for every
$s\in\mathcal{S}=(m+n\mathbf{Z})\cup (m+1+n\mathbf{Z})$, making
commute the following diagrams:

\vspace*{0.5cm}

\setlength{\unitlength}{1mm}
\begin{picture}(140,30)
\put(10,25){$...I^{m+(k-1)n+1}$} \put(32,26){\vector(1,0){18}}
\put(52,25){$I^{m+kn}$} \put(62,26){\vector(1,0){18}}
\put(82,25){$I^{m+kn+1}...$}

\put(10,4){$...J^{m+(k-1)n+1}$} \put(32,5){\vector(1,0){18}}
\put(52,4){$J^{m+kn}$} \put(62,5){\vector(1,0){18}}
\put(82,4){$J^{m+kn+1}...$}

\put(14,23){\vector(0,-1){12}} \put(53,23){\vector(0,-1){12}}
\put(83,23){\vector(0,-1){12}}

\put(39,28){$d^{n-1}$} \put(71,28){$d$}  \put(39,2){$d^{n-1}$}
\put(71,2){$d$}

\put(11,16){$f$} \put(50,16){$f$} \put(80,16){$f$}

\end{picture}

where $d^{n-1}$ means the composition of $n-1$ consecutive
'differentials' of $I^\cdot$.  By \cite{MVS}[Lemma 2.8], that
gives uniquely determined morphisms in $Mod_{\Lambda_0}$
$g_s:M_s\longrightarrow N_s$, for every $s\in (m+n\mathbf{Z})\cup
(m+1+n\mathbf{Z})$, such that $f^s=\mathcal{H}om (A,g_s)[s]$. The
commutativity of the above diagram translates into the
commutativity of the following diagrams, for every
$k\in\mathbf{Z}$:

\vspace*{0.5cm}

\setlength{\unitlength}{1mm}
\begin{picture}(140,30)

\put(15,25){$Hom_{\Lambda_0}(\Lambda_{n-1},M_{m+(k-1)n+1})$}
\put(60,26){\vector(1,0){10}} \put(72,25){$M_{m+kn}$}

\put(15,4){$Hom_{\Lambda_0}(\Lambda_{n-1},N_{m+(k-1)n+1})$}
\put(60,5){\vector(1,0){10}} \put(72,4){$N_{m+kn}$}

\put(35,23){\vector(0,-1){16}} \put(75,23){\vector(0,-1){16}}
\put(30,15){$g_*$}  \put(72,15){$g$}

\put(62,27){$d^{n-1}$} \put(62,2){$d^{n-1}$}

\end{picture}

\vspace*{0.5cm}

\setlength{\unitlength}{1mm}
\begin{picture}(140,30)

\put(15,25){$Hom_{\Lambda_0}(\Lambda_{1},M_{m+kn})$}
\put(50,26){\vector(1,0){20}} \put(72,25){$M_{m+kn+1}$}

\put(15,4){$Hom_{\Lambda_0}(\Lambda_{1},M_{m+kn})$}
\put(50,5){\vector(1,0){20}} \put(72,4){$M_{m+kn+1}$}

\put(35,23){\vector(0,-1){16}} \put(75,23){\vector(0,-1){16}}
\put(30,15){$g_*$}  \put(72,15){$g$}

\put(60,27){$d$} \put(60,2){$d$}

\end{picture}

where $d^{n-1}(h)=\sum_{p\in Q_{n-1}}h(p)p^o$ and
$d(h)=\sum_{\alpha\in Q_1}h(\alpha )\alpha^o$, respectively. Now
we consider the obvious adaptation of \cite{MVS}[Lemma 2.9], which
is also true replacing $Q_1$ by $Q_{n-1}$ due to the fact that
$\Lambda^!_i=KQ_i^{op}$ for $0\leq i<n$. Then the commutativity of
the last two diagrams is equivalent to say that $g(xa)=g(x)a$
whenever $a\in \Lambda^!_1\cup \Lambda^!_{n-1}$ and both $deg(x)$
and $deg(xa)$ belong to $\mathcal{S}=(m+n\mathbf{Z})\cup
(m+1+n\mathbf{Z})$. If now $a\in \Lambda^!_n$ and $x\in M_{m+kn}$,
then we have $a=\sum_i\alpha_i\beta_i$, where $\alpha_i\in
\Lambda^!_1=KQ_1^{op}$ and $\beta_i\in
\Lambda^!_{n-1}=KQ_{n-1}^{op}$. Then $g(xa)=\sum_i
g(x(\alpha_i\beta_i))=\sum_i g((x\alpha_i)\beta_i)$. Now
$deg(x\alpha_i)=m+kn+1\in\mathcal{S}$ and the upper diagram above
gives $\sum_i g((x\alpha_i)\beta_i)=\sum_i g(x\alpha_i)\beta_i$,
while the lower diagram gives that the latter expression equals
$\sum_i (g(x)\alpha_i)\beta_i)=\sum_i
g(x)(\alpha_i\beta_i)=g(x)\sum_i\alpha_i\beta_i=g(x)a$, using the
associative property of the $A$-module $N$. One can proceed in an
analogous way when $x\in M_{m+kn+1}$ and $a\in \Lambda^!_n$, but
taking a decomposition $a=\sum_i\beta_i\gamma_i$, where
$\beta_i\in \Lambda^!_{n-1}=KQ_{n-1}^{op}$ and $\gamma_i\in
\Lambda^!_1=KQ_1^{op}$. According to Remark \ref{remark}, we
conclude that $g:M_\mathcal{S}\longrightarrow N_\mathcal{S}$ is a
morphism in $Gr_{A_\mathcal{U}}$.
 Since, by definition, both
$M_{\mathcal{S}}$ and $N_{\mathcal{S}}$ are liftable with respect
to $(-)_{\mathcal{S}}$ and generated in degrees belonging to
$(\mathcal{S}:\mathcal{U})$, Theorem 2.7 of \cite{MVS2} tells us
that there exists a unique morphism $\eta :M\longrightarrow N$ in
$Gr_{\Lambda^!}$ such that $\eta_{\mathcal{S}}=g$. It is now a
mere routine to check that $\bar{f}=\upsilon (\eta
):I^\cdot\longrightarrow J^\cdot$ is a morphism in
$_n\mathcal{L}\mathcal{C}^*_\Lambda$ such that $H(\bar{f})=f$.

\end{proof}

We can now put together all the pieces of the puzzle. Recall from
\cite{MVS2} that if $X\in Gr_{\Lambda^!_\mathcal{U}}$, with
$Supp(X)\subseteq\mathcal{S}$, and  $\mu_{m+kn,1}:X_{m+kn}\otimes
\Lambda^!_1\longrightarrow X_{m+kn+1}$ and
$\mu_{m+kn,1n}:X_{m+kn}\otimes \Lambda^!_n\longrightarrow
X_{m+(k+1)n}$ are the multiplication maps, then, working within
the right $\Lambda^!$-module $X_{m+n\mathbf{Z}}\otimes\Lambda^!$,
it makes sense to consider $Ker(\mu_{m+kn,1})\Lambda^!_{n-1}$,
which is a $\Lambda_0$-submodule of $X_{m+kn}\otimes\Lambda^!_n$.
We now have:

\begin{teor} \label{grmodules-to-2complexes}
Let $n\geq 2$ be a positive integer and  consider the subsets
$\mathcal{U}=n\mathbf{Z}\cup (n\mathbf{Z}+1)$ and
$\mathcal{S}=m+\mathcal{U}$, where $m\in\mathbf{Z}$. Let us assume
that $\Lambda $ is a graded factor of a path algebra with
relations of degree $\geq n$. There is an equivalence between the
following categories:

\begin{enumerate}
\item The full subcategory $\mathcal{L}(\mathcal{S},\mathcal{U})$
of $Gr_{\Lambda^!_\mathcal{U}}$ with objects those $X$ which are
generated in degrees belonging to $(\mathcal{S}:\mathcal{U})$ and
satisfying that $Ker(\mu_{m+kn,1})\Lambda^!_{n-1}\subseteq
Ker(\mu_{m+kn,n})$, for all $k\in\mathbf{Z}$ \item The full
subcategory $\mathcal{Y}(\mathcal{S},\mathcal{U})$ of
$C(Gr_\Lambda )$ whose objects are the $H_m$-liftable
(2-)complexes $(\tilde{I}^\cdot ,\tilde{d}^\cdot )$ satisfying
contition a) below, in case $n=2$, and both conditions a) and b),
in case $n>2$:

\begin{enumerate}
\item $\tilde{I}^j$ is almost injective and cogenerated in degree
$-\delta_m(j)$, for every $j\in\mathbf{Z}$ \item
$Soc^{gr}(\tilde{I}^{2j+1})\subseteq Im(\tilde{d}^{2j})$, for all
$j\in\mathbf{Z}$
\end{enumerate}
\end{enumerate}
Moreover, $\mathcal{L}(\mathcal{S},\mathcal{U})$ contains all the
graded $A_\mathcal{U}$-modules presented in degrees belonging to
$(\mathcal{S}:\mathcal{U})$.
\end{teor}

\begin{proof}
In case $n=2$ everything is trivial, and is actually a shifting
version of the equivalence
$Gr_{\Lambda^!}\stackrel{\cong}{\longrightarrow}\mathcal{LC}^*_\Lambda$
of \cite{MVS}[Theorem 2.10]. Indeed instead of the functor $\nu
=\nu_\Lambda$ there given, ours here is the composition
$Gr_{\Lambda^!}\stackrel{\nu}{\longrightarrow}C(Gr_\Lambda
)\stackrel{?[m]}{\longrightarrow}C(Gr_\Lambda )$.

For the case $n>2$, notice that
$(\mathcal{S}:\mathcal{U})=m+n\mathbf{Z}$. From Proposition
\ref{transfer-of-torsion-and-quotient-categories}, Proposition
\ref{from-ncomplex-to-2complex1} and \cite{MVS2}[Theorem 2.7] we
get the following diagram, where all the thick arrows are
equivalences of categories:

\vspace*{0.5cm}

\setlength{\unitlength}{1mm}
\begin{picture}(140,30)

\put(30,25){$\mathcal{G}(\mathcal{S},\mathcal{U})$}
\put(45,26){\vector(1,0){15}}
\put(62,25){$\mathcal{G}^*(\mathcal{S},\mathcal{U})$}

\put(30,4){$\mathcal{L}(\mathcal{S},\mathcal{U})$}
  \put(62,4){$\mathcal{Y}(\mathcal{S},\mathcal{U})$}

\put(51,28){$\upsilon_\Lambda$}

\put(36,23){\vector(0,-1){15}}  \put(66,23){\vector(0,-1){15}}
\put(28,15){$(-)_\mathcal{S}$}  \put(68,15){$H$}

 \put(50,5){$\rightsquigarrow$}

 \put(37,15){$\cong$}   \put(63,15){$\cong$}   \put(51,23){$\cong$}

\end{picture}

 Then the squig arrow making commute the diagram is also an equivalence.
\end{proof}

\begin{rem} \label{remark-important}

1) We  make explicit, in case $n>2$,  the definition of the
equivalence
$\mathcal{L}(\mathcal{S},\mathcal{U})\stackrel{\cong}{\longrightarrow}\mathcal{Y}(\mathcal{S},\mathcal{U})$.
Assume that $\Lambda =\frac{KQ}{I}$, where $I$ is generated by
$\rho\subseteq \coprod_{i\geq n}KQ_i$. Then
$\Lambda^!=\frac{KQ^{op}}{<I_n^\perp
>}$. If $X\in \mathcal{L}(\mathcal{S},\mathcal{U})$, then we have
a morphism of $\Lambda_0$-modules $\xi =\xi_k:X_{m+kn+1}\otimes
\Lambda^!_{n-1}=X_{m+kn+1}\otimes KQ_{n-1}^{op}\longrightarrow
X_{m+(k+1)n}$ defined as follows. Any element of
$X_{m+kn+1}\otimes KQ_{n-1}^{op}$ can be written as $\sum_{p\in
Q_{n-1}}x_p\otimes p^o$, for uniquely determined $x_p\in
X_{m+kn+1}e_{t(p)}$, where $e_i$ denotes the idempotent of $KQ_0$
corresponding to the vertex $i\in Q_0$. Since
$X_{m+kn+1}=X_{m+kn}\Lambda^!_1=X_{m+kn}KQ_1^{op}$, we can write
$x_p$ as $x_p=\sum_{\alpha\in Q_1, o(\alpha )=t(p)}x_{\alpha
,p}\alpha^o$. Then we put $\xi (\sum_{p\in Q_{n-1}}x_p\otimes
p^o)=\sum_{\alpha ,p}x_{\alpha ,p}(\alpha^op^o)=\sum_{\alpha
,p}x_{\alpha ,p}(p\alpha )^o$. That $\xi$ is well-defined follows
from the inclusions $Ker(\mu_{m+kn,1})\Lambda^!_{n-1}\subseteq
Ker(\mu_{m+kn,n})$ given by the theorem. We put then
$\tilde{I}^{2j}=\mathcal{H}om_{\Lambda_0}(\Lambda
,X_{m+jn})[m+jn]$ and
$\tilde{I}^{2j+1}=\mathcal{H}om_{\Lambda_0}(\Lambda
,X_{m+jn+1})[m+jn+1]$, for all $j\in\mathbf{Z}$. The differential
$\tilde{d}^{2j}$ is completely identified by the map
$Hom_{\Lambda_0}(\Lambda_1,X_{m+jn})\longrightarrow X_{m+jn+1}$
which takes $f\rightsquigarrow\sum_{\alpha\in Q_1}f(\alpha
)\alpha^o$ and the differential
$\tilde{d}^{2j+1}:\tilde{I}^{2j+1}={H}om_{\Lambda_0}(\Lambda
,X_{m+jn+1})[m+jn+1]\longrightarrow {H}om_{\Lambda_0}(\Lambda
,X_{m+(j+1)n})[m+(j+1)n]=\tilde{I}^{2j+2}$ is completely
identified by the map
$Hom_{\Lambda_0}(\Lambda_{n-1},X_{m+jn+1})=Hom_{\Lambda_0}(KQ_{n-1},X_{m+jn+1})
\longrightarrow X_{m+(j+1)n}$ which takes $f\rightsquigarrow
\sum_{p\in Q_{n-1}}\xi (f(p)\otimes p^o)$. It is a mere routine to
check that the equivalence given by last theorem takes $X$ onto
the here defined $(\tilde{I}^\cdot ,\tilde{d}^\cdot)$.

2) If, instead of taking the functor $\nu=\nu_\Lambda$ of
Proposition \ref{gradedmodules-ncomplexes2} as basis of our
arguments, one takes the functor $\Psi =\Psi_\Lambda$ of
Proposition \ref{gradedmodules-ncomplexes1}, then one gets results
dual to those in sections 3 and 4. In the case of locally finite
graded modules,  we can  alternatively see that by using the
canonical duality $D$. Indeed, we have $D(N\otimes\Lambda
)\cong\mathcal{H}om_{\Lambda_0}(\Lambda ,N)$, for all $N\in
mod_{\Lambda_0}$, and then the commutativity of the following
diagram follows, where $\upsilon_\Lambda$ is taken for graded left
$\Lambda^!$-modules:

\vspace*{0.5cm}

\setlength{\unitlength}{1mm}
\begin{picture}(140,30)

\put(30,25){$_{\Lambda^!}lfgr$} \put(45,26){\vector(1,0){15}}
\put(62,25){$_n\mathcal{L}c^*_{\Lambda^{op}}$}

\put(30,4){$lfgr_{\Lambda^!}$} \put(45,5){\vector(1,0){15}}
  \put(62,4){$_n\mathcal{L}c_\Lambda$}

\put(51,28){$\upsilon_\Lambda$} \put(51,2){$\Psi_\Lambda$}

\put(36,23){\vector(0,-1){15}}  \put(66,23){\vector(0,-1){15}}
\put(31,15){$D$}  \put(68,15){$D$}

\put(37,15){$\cong^o$}  \put(62,15){$\cong^o$}
\put(51,23){$\cong$} \put(51,6){$\cong$}

\end{picture}

 We make
explicit, without proof,  the dual of Theorem
\ref{grmodules-to-2complexes}, leaving the rest for the reader.

\end{rem}

We have a canonical additive functor
$G=G_m:_n\mathcal{L}C_\Lambda\longrightarrow C(Gr_\Lambda )$
defined as follows. If $(P^\cdot ,d^\cdot )\in
_n\mathcal{L}C_\Lambda$ then its image by $G$ is $(\tilde{P}^\cdot
,\tilde{d}^\cdot )$, where $\tilde{P}^j=P^{-\delta_m (-j)}$ and
the differentials are
$\tilde{d}^{2j-1}=d^{-m+jn-1}:\tilde{P}^{2j-1}=P^{-m+jn-1}\longrightarrow
P^{-m+jn}=\tilde{P}^{2j}$ and $\tilde{d}^{2j}=d^{-m+(j+1)n-2}\circ
...\circ d^{-m+nj}:\tilde{P}^{2j}=P^{-m+jn}\longrightarrow
P^{-m+(j+1)n-1}=\tilde{P}^{2j+1}$, for all $j\in\mathbf{Z}$. To
state the desired dual, for every $X\in
Gr_{\Lambda^!_\mathcal{U}}$ with $Supp(X)\subseteq
-\mathcal{S}=(-m+n\mathbf{Z})\cup (-m-1+n\mathbf{Z})$, we consider
the comultiplications $\Delta_{s,u}:X_{-s-u}\longrightarrow
X_{-s}\otimes KQ_u$, for all $(s,u)\in(\mathcal{S},\mathcal{U})$,
such that $s+u\in\mathcal{S}$ and $u\geq 0$. By definition, one
has  $\Delta_{s,u}(x)=\sum_{p\in Q_u}x\bar{p}^o\otimes p$, where
$\bar{p}^o$ is the image of $p^o$ by the canonical projection
$KQ_u^{op}\twoheadrightarrow\Lambda^!_u$. We are now ready to
state the dual of last theorem.

\begin{teor} \label{grmodules-to-2complexes2}

Let $n\geq 2$ be a positive integer and and consider the subsets
$\mathcal{U}=n\mathbf{Z}\cup (n\mathbf{Z}+1)$ and
$\mathcal{S}=m+\mathcal{U}$, where $m\in\mathbf{Z}$. Let us assume
that $\Lambda$ is a graded factor of a path algebra with relations
of degree $\geq n$. There is an equivalence between the following
categories:

\begin{enumerate}
\item The full subcategory
$\mathcal{L}^o(\mathcal{S},\mathcal{U})$ of
$Gr_{\Lambda^!_\mathcal{U}}$ with objects those $X$ which are
cogenerated in degrees belonging to $-(\mathcal{S}:\mathcal{U})$
and satisfy that the following diagram in $Mod_{\Lambda_0}$ can be
completed for all $k\in\mathbf{Z}$:

\vspace*{0.5cm}

\setlength{\unitlength}{1mm}
\begin{picture}(140,30)

\put(15,25){$X_{-m-(k+1)n}$} \put(35,26){\vector(1,0){30}}
\put(68,25){$X_{-m-kn}\otimes KQ_n$}

\put(9,4){$X_{-m-kn-1}\otimes KQ_{n-1}$}
\put(42,5){\vector(1,0){23}} \put(68,4){$X_{-m-kn}\otimes
KQ_1\otimes KQ_{n-1}$}

\put(80,23){\vector(0,-1){16}}

\put(50,28){$\Delta$} \put(47,2){$\Delta\otimes 1$}
\put(82,15){$\cong$}

\put(24,15){$\downarrow$}

\end{picture}

\item The full subcategory $\mathit{Y}^o(\mathcal{S},\mathcal{U})$
of $C(Gr_\Lambda )$ whose objects are the $G_m$-liftable
(2-)complexes $(\tilde{P}^\cdot ,\tilde{d}^\cdot )$ satisfying
condition a) below, in case $n=2$, and conditions a) and b), in
case $n>2$:

\begin{enumerate}
\item $\tilde{P}^j$ is a  projective graded $\Lambda$-module
generated in degree $\delta_m (-j)$, for every $j\in\mathbf{Z}$
\item $(Ker\tilde{d}^{2k-1})\subseteq\tilde{P}^{2k-1}\cdot
J^{gr}(\Lambda )$, for all $k\in\mathbf{Z}$
\end{enumerate}
\end{enumerate}
\end{teor}

In case $\Lambda$ is a $n$-Koszul algebra with Yoneda algebra
$E=E(\Lambda )$ we know from \cite{GMMVZ}[Theorem 9.1] (see also
\cite{BM}[Proposition 3.1]) that there is an algebra isomorphism
$\varphi
:E\stackrel{\cong}{\longrightarrow}\Lambda^!_\mathcal{U}$, which
we fix from now on,  such that $\varphi (E_j)=\Lambda^!_{\delta
(j)}$, for all $j\in\mathbf{Z}$. We see it as an identification
and, abusing of notation, we write $E_j=\Lambda^!_{\delta (j)}$,
for all $j\in\mathbf{Z}$. Hence, if $V\in Gr_E$, we have
multiplication maps $\tilde{\mu}_{2k,1}:V_{2k}\otimes
E_1=V_{2k}\otimes\Lambda^!_1\longrightarrow V_{2k+1}$ and
$\tilde{\mu}_{2k,2}:V_{2k}\otimes
E_2=V_{2k}\otimes\Lambda^!_n\longrightarrow V_{2k+2}$ and we can
take $Ker(\tilde{\mu}_{2k,1})\Lambda^!_{n-1}$, which is an
$\Lambda_0-\Lambda_0-$sub-bimodule of
$V_{2k}\otimes\Lambda^!_n=V_{2k}\otimes E_2$. We then have the
following  consequence of Theorem \ref{grmodules-to-2complexes}:

\begin{cor} \label{modules-over-Yoneda}
Let $n\geq 2$ be a positive integer and consider the subsets
$\mathcal{U}=n\mathbf{Z}\cup (n\mathbf{Z}+1)$ and
$\mathcal{S}=m+\mathcal{U}$, where $m\in\mathbf{Z}$. If $\Lambda
=\oplus_{i\geq 0}\Lambda_i$ is a $n$-Koszul algebra and $E$ is its
Yoneda algebra, then there is an equivalence between:

\begin{enumerate}

\item  The full subcategory $\mathcal{L}_E$   of $Gr_E$, which, in
case $n=2$,  coincides with $Gr_E$  and, in case $n>2$,  has as
objects those $V\in Gr_E$ which are generated in even degrees and
satisfy that $Ker(\tilde{\mu}_{2k,1})\Lambda^!_{n-1}\subseteq
Ker(\tilde{\mu}_{2k,2})$, for all $k\in\mathbf{Z}$ \item The full
subcategory $\mathcal{Y}(\mathcal{S},\mathcal{U})$ of
$C(Gr_\Lambda )$ defined in Theorem \ref{grmodules-to-2complexes}
\end{enumerate}

Moreover, $\mathcal{L}_E$ contains all the graded $E$-modules
presented in even degrees.
\end{cor}

\begin{proof}

In case $n=2$, one has $E=\Lambda^!$,  $\mathcal{L}_E=Gr_E$ and
$\mathcal{Y}(\mathcal{S},\mathcal{U})$ is the full subcategory of
$C(Gr_\Lambda )$ consisting of those complexes $(\tilde{I}^\cdot
,\tilde{d}^\cdot )$ such that $\tilde{I}^j$ is almost injective
and cogenerated in degree $-m-j$, for every $j\in\mathbf{Z}$. The
desired equivalence of categories is then the composition of the
equivalences $\nu_\Lambda
:Gr_E=Gr_{\Lambda^!}\stackrel{\cong}{\longrightarrow}\mathcal{LC}^*_\Lambda$
and
$?[m]:\mathcal{LC}^*_\Lambda\stackrel{\cong}{\longrightarrow}\mathcal{Y}(\mathcal{S},\mathcal{U})$.

We next assume $n>2$. Then, from \cite{MVS2}[Theorem 2.7,
Corollary 4.10 and Remark 4.11] we get equivalences of categories
$\mathcal{L}_E\cong\mathcal{G}(\mathcal{S},\mathcal{U})
\cong\mathcal{L}(\mathcal{S},\mathcal{U})$ and, using Theorem
\ref{grmodules-to-2complexes}, also an equivalence
$\mathcal{L}(\mathcal{S},\mathcal{U})\cong\mathcal{Y}(\mathcal{S},\mathcal{U})$.

\end{proof}

\begin{ejem}
Let us take a positive integer $n>2$  and consider the truncated
algebra $\Lambda =KQ/<Q_n>$, so that $\Lambda$ is $n$-Koszul and
$\Lambda^!=KQ^{op}$. We give the equivalence of last corollary
when  $m=0$ (i.e. $\mathcal{S}=\mathcal{U}=n\mathbf{Z}\cup
(n\mathbf{Z}+1)$). Then, as ungraded algebras,
$E=\Lambda^!_\mathcal{U}=(KQ^{op})_\mathcal{U}$.  The classical
grading of $E$ is given by assigning degree 1 to the $\alpha^o\in
Q_1^o$ and degree 2 to the $p^o\in Q_n^o$. Then the subcategory
$\mathcal{L}_E$ consists of those $V\in Gr_E$ which are generated
in even degrees and satisfy that if $\sum_{\alpha\in
Q_1}x_\alpha\alpha^o =0$, for a family of elements $(x_\alpha )$
in $V_{2j}$, then $\sum_{\alpha\in Q_1}x_\alpha(\alpha^oq^o) =0$
for all $q\in Q_{n-1}$. Using Remark \ref{remark-important}(1)we
then get that the equivalence
$\mathcal{L}_E\stackrel{\cong}{\longrightarrow}\mathcal{Y}(\mathcal{U},\mathcal{U})$
maps $V$ onto the cochain complex $(\tilde{I}^\cdot
,\tilde{d}^\cdot )$ defined as follows: i)
$\tilde{I}^{2j}=\mathcal{H}om_{\Lambda_0}(\Lambda ,V_{2j})[jn]$
and $\tilde{I}^{2j+1}=\mathcal{H}om_{\Lambda_0}(\Lambda
,V_{2j+1})[jn+1]$, for all $j\in\mathbf{Z}$; ii) the differentials
$\tilde{I}^{2j}=\mathcal{H}om_{\Lambda_0}(\Lambda
,V_{2j})[jn]\longrightarrow\mathcal{H}om_{\Lambda_0}(\Lambda
,V_{2j+1})[jn+1]=\tilde{I}^{2j+1}$ and
$\tilde{I}^{2j+1}=\mathcal{H}om_{\Lambda_0}(\Lambda
,V_{2j+1})[jn+1]\longrightarrow \mathcal{H}om_{\Lambda_0}(\Lambda
,V_{2j+2})[(j+1)n]=\tilde{I}^{2j+2}$ are defined by the formulas
$\tilde{d}(f)(a)=\sum_{\alpha\in Q_1}f(a\alpha )\alpha^o$ and
$\tilde{d}(g)(a)=\sum_{p\in Q_{n-1}}g(ap)\cdot p^o$, for all
$a\in\Lambda_i$, $f\in Hom_{\Lambda_0}(\Lambda_{i+1},V_{2j})$ and
$g\in Hom_{\Lambda_0}(\Lambda_{i+n-1},V_{2j+1})$. The last
multiplication  $\cdot
:V_{2j+1}\otimes\Lambda^!_{n-1}\longrightarrow V_{2j+2}$ is given
by $(\sum_{\alpha\in Q_1}x_\alpha\alpha^o)\cdot
p^o=\sum_{\alpha\in Q_1}x_\alpha (\alpha^op^o)=\sum_{\alpha\in
Q_1}x_\alpha (p\alpha)^o$, for all $p\in Q_{n-1}$ and all
$x=\sum_{\alpha\in Q_1}x_\alpha\alpha^o\in X_{2j+1}$, where the
$x_\alpha$ belong to $X_{2j}$.
\end{ejem}

If  $M\in lfgr_\Lambda$ we shall say that $M$ is n-{\bf coKoszul}
when it is cogenerated in degree $0$ and  its minimal (almost)
injective graded resolution $\tilde{I}^\cdot _M$ satisfies that
$\tilde{I}^j$ is cogenerated in degree $-\delta (j)$, for every
$j\geq 0$, where $\delta =\delta_0$. In particular,
$\tilde{I}^\cdot _M$ satisfies conditions a) and b) of Theorem
\ref{grmodules-to-2complexes} for $m=0$. It is easy to see that
the assigment $M\rightsquigarrow\tilde{I}^\cdot _M$ yields a fully
faithful embedding of the category $\mathcal{K}_n^o(\Lambda )$ of
n-coKoszul $\Lambda$-modules into $C(lfgr_\Lambda)$. We denote by
$\mathcal{K}_n(\Lambda)$ the full subcategory of $_\Lambda lfgr$
formed by the $n$-Koszul modules. It would be interesting to have
an answer to the following question:

\begin{quesrem}
Let $\Lambda$ be a $n$-Koszul algebra ($n>2$). Which are the
$n$-co-Koszul modules $M$ such that $\tilde{I}^\cdot_M$ is
$H_0$-liftable (and, hence, belongs to
$\mathcal{Y}(\mathcal{U},\mathcal{U})$)?.  If we denote the
corresponding full subcategory of $lfgr_\Lambda$ by
$\mathcal{K}^{o}_H(\Lambda )$, then its image by the canonical
duality $D:lfgr_\Lambda\longrightarrow _\Lambda lfgr$, which we
denote $\mathcal{K}_G(\Lambda )$,  consists of those (locally
finite) $n$-Koszul modules whose minimal graded projective
resolution is $G_0$-liftable (see Theorem
\ref{grmodules-to-2complexes2}). Now, going backward in Corollary
\ref{modules-over-Yoneda}, we get a contravariant fully faithful
embedding $\mathcal{K}_G(\Lambda
)\stackrel{\cong^o}{\longrightarrow}\mathcal{K}^{o}_H(\Lambda
)\rightarrowtail\mathcal{L}_E\hookrightarrow Gr_E$ which takes
indecomposable projective graded $\Lambda$-modules  onto simple
graded $E$-modules  and simple graded $\Lambda$-modules onto
indecomposable projective graded $E$-modules. The last assertion
follows from \cite{BM}[Section 3], where the authors prove that
the minimal projective resolution of $\Lambda_0$ in $_\Lambda Gr$
is obtained by contraction of the Koszul n-complex, and is thereby
$G$-liftable.

\end{quesrem}


\begin{thebibliography}{9999}



\bibitem{BGS} {\sc BEILINSON, A.; GINZBURG, V.; SOERGEL, W.}:
Koszul duality patterns in representation theory. J. Amer. Math.
Soc. {\bf 9}(2) (1996), 473-526

\bibitem{B}{\sc BERGER, R.}: Koszulity for nonquadratic algebras.
J. Algebra {\bf 239} (2001), 705-734

\bibitem{BDVW} {\sc BERGER, R.; DUBOIS-VIOLETTE, M.; WAMBST, M.}:
Homogeneous algebras. J. Algebra {\bf 261} (2003), 172-185

\bibitem{BM} {\sc BERGER, R.; MARCONNET, N.}: Koszul and
Gorenstein properties for homogeneous algebras. Preprint,   arXiv:
math.QA/0310070 v2

\bibitem{BG} {\sc BERNSTEIN, J.; GELFAND, S.}: Algebraic vector
bundles and projective spaces. Appendix to Russian translation of
M. Schneider, 'Holomorphic vector bundles on $\mathbf{P}^n$'. Sem.
Bourbaki {\bf 530} (1980), 80-102


\bibitem{GM} {\sc GELFAND, S.I.; MANIN, Y.I.}: "Methods of
homological algebra". Springer-Verlag (1996)

\bibitem{GMMVZ} {\sc GREEN, E.L.; MARCOS, E.N.; MARTINEZ-VILLA, R.; ZHANG,
P.}: D-Koszul algebras. J. Pure and Appl. Algebra {\bf 193}
(2004), 359-378


\bibitem{K} {\sc KELLER, B.}: Introduction to
$\mathbf{A}_\propto$-algebras and modules. Homology, Homotopy and
Appl. {\bf 3} (2001), 1-35

\bibitem{MVS}{\sc MARTINEZ VILLA, R.; SAORIN, M.}: Koszul
equivalences and dualities. Pacific J. Math. {\bf 204}(2) (2004),
359-378

\bibitem{MVS2}{\sc MARTINEZ VILLA, R.; SAORIN, M.}: Killing of
supports on graded algebras. Preprint

\bibitem{NVO} {\sc NASTASESCU, C.; VAN OYSTAEYEN, F.}: "Graded
Ring Theory". North-Holland (1982)



\end{thebibliography}
\end{document}